\documentclass{amsart}

\theoremstyle{definition}

\theoremstyle{remark}

\numberwithin{equation}{section}

\begin{document}

\title{ON THE PROXIMINALITY OF RIDGE FUNCTIONS}

\author{Vugar E. Ismailov}

\address{Mathematics and Mechanics Institute,
Azerbaijan National Academy of Sciences, Az-1141, Baku,
Azerbaijan}

\email{vugaris@mail.ru}

\thanks{This research was supported by INTAS under Grant
06-1000015-6283}

\subjclass[2000]{41A30, 41A50, 41A63}

\keywords{Ridge function; Extremal element; Proximinality; Path;
Orbit}

\begin{abstract}
Using two results of Garkavi, Medvedev and Khavinson [7], we give
sufficient conditions for proximinality of sums of two ridge
functions with bounded and continuous summands in the spaces of
bounded and continuous multivariate functions respectively. In the
first case, we give an example which shows that the corresponding
sufficient condition cannot be made weaker for some subsets of
$\mathbb{R}^{n}$. In the second case, we obtain also a necessary
condition for proximinality. All the results are furnished with
plenty of examples. The results, examples and following
discussions naturally lead us to a conjecture on the proximinality
of the considered class of ridge functions. The main purpose of
the paper is to draw readers' attention to this conjecture.
\end{abstract}

\maketitle

\section*{0. Introduction}

In multivariate approximation theory, special functions called
\textit{ridge functions} are widely used. A ridge function is a
multivariate function of the form $g\left( \mathbf{a}\cdot
\mathbf{x}\right) $, where $g$ is a univariate function,
$\mathbf{a}=\left( a_{1},...,a_{n}\right) $ is a vector
(direction) different from zero, $\mathbf{x}=\left(
x_{1},...,x_{n}\right) $ is the variable and $\mathbf{a}\cdot
\mathbf{x}$ is the inner product. In other words, a ridge function
is a composition of a univariate function with a linear functional
over $\mathbb{R}^{n}.$ These functions arise naturally in various
fields. They arise in partial differential equations (where they
are called \textit{plane waves} [15]), in computerized tomography
(see, e.g., [19,22]; the name ridge function was coined by Logan
and Shepp[19] in one of the seminal papers on tomography), in
statistics (especially, in the theory of projection pursuit and
projection regression; see, e.g., [4,11]). Ridge functions are
also the underpinnings of many central models in neural networks
which has become increasing more popular in computer science,
statistics, engineering, physics, etc. (see [24] and references
therein). We refer the reader to Pinkus [23] for various
motivations for the study of ridge functions and ridge function
approximation.

Let $E$ be a normed linear space and $F$ be its subspace. We say
that $F$ is proximinal in $E$ if for any element $e\in E$ there
exists at least one element $f_{0}\in F$ such that

\begin{equation*}
\left\Vert e-f_{0}\right\Vert =\inf_{f\in F}\left\Vert
e-f\right\Vert .
\end{equation*}

In this case, the element $f_{0}$ is said to be extremal to $e$.

Although at present there are a great deal of interesting papers
devoted to the approximation by ridge functions (see, e.g.,
[2,3,5,9,12,13,17,18,20,24,25]), some problems of this
approximation have not been solved completely yet. In the
following, we are going to deal with one of such problems, namely
with the problem of proximinality of the set of linear
combinations of ridge functions in the spaces of bounded and
continuous functions respectively. This problem will be considered
in the simplest case when the class of approximating functions is
the set

\begin{equation*}
\mathcal{R}=\mathcal{R}\left(
\mathbf{a}^{1}{,}\mathbf{a}^{2}\right) ={ \left\{ {g_{1}\left(
\mathbf{a}^{1}{\cdot }\mathbf{x}\right) +g_{2}\left(
\mathbf{a}^{2}{\cdot }\mathbf{x}\right)
:g}_{i}:{{\mathbb{R\rightarrow R}} ,i=1,2}\right\} }.
\end{equation*}
Here $\mathbf{a}^{1}${and }$\mathbf{a}^{2}$ are fixed directions
and we vary over ${g}_{i}$. It is clear that this is a linear
space. Consider the following three subspaces of $\mathcal{R}$.
The first is obtained by taking only bounded sums ${g_{1}\left(
\mathbf{a}^{1}{\cdot }\mathbf{x}\right) +g_{2}\left(
\mathbf{a}^{2}{\cdot }\mathbf{x}\right) }$ over some set $X$ in
$\mathbb{R}^{n}.$ We denote this subspace by $\mathcal{R}_{a}(X)$.
The second and the third are subspaces of $\mathcal{R}$ with
bounded and continuous summands $g_{i}\left( \mathbf{a}^{i}\cdot
\mathbf{x}\right) ,~i=1,2,$ on $X$ respectively. These subspaces
will be denoted by $\mathcal{R }_{b}(X)$ and $\mathcal{R}_{c}(X).$
In the case of $\mathcal{R}_{c}(X),$ the set $X$ is considered to
be compact.

Let $B(X)$ and $C(X)$ be the spaces of bounded and continuous
multivariate functions over $X$ respectively. What conditions must
one impose on $X$ in order that the sets $\mathcal{R}_{a}(X)$ and
$\mathcal{R} _{b}(X)$ be proximinal in $B(X)$ and the set
$\mathcal{R}_{c}(X)$ be proximinal in $C(X)$? We are also
interested in necessary conditions for proximinality. It follows
from one result of Garkavi, Medvedev and Khavinson (see theorem1
[7]) that $\mathcal{R}_{a}(X)$ is proximinal in $B(X)$ for all
subsets $X$ of $\mathbb{R}^{n}$. There is also an answer (see
theorem 2 [7]) for proximinality of $\mathcal{R}_{b}(X)$ in
$B(X)$. This will be discussed in Section 1. Is the set
$\mathcal{R}_{b}(X)$ always proximinal in $B(X)$? There is an an
example of a set $X\subset \mathbb{R}^{n}$ and a bounded function
$f$ on $X$ for which there does not exist an extremal element in $
\mathcal{R}_{b}(X)$.

In Section 2, we will obtain sufficient conditions for the
existence of extremal elements from $\mathcal{R}_{c}(X)$ to an
arbitrary function $f$ $\in$ $C(X)$. Based on one result of
Marshall and O'Farrell [21], we will also give a necessary
condition for proximinality of $\mathcal{R}_{c}(X) $ in $C(X)$.
All the theorems, following discussions and examples of the paper
will lead us naturally to a conjecture on the proximinality of the
subspaces $\mathcal{R}_{b}(X)$ and $\mathcal{R}_{c}(X)$ in the
spaces $B(X)$ and $C(X)$ respectively.

At the end of this section, we want to draw the readers attention
to the more general case in which the number of directions is more
than two. In this case, the set of approximating functions is

\begin{equation*}
\mathcal{R}\left( \mathbf{a}^{1},...,\mathbf{a}^{r}\right)
=\left\{ \sum\limits_{i=1}^{r}g_{i}\left( \mathbf{a}^{i}\cdot
\mathbf{x}\right) :g_{i}:\mathbb{R}\rightarrow
\mathbb{R},i=1,...,r\right\} .
\end{equation*}

In a similar way as above, one can define the sets $\mathcal{R}
_{a}(X)$, $\mathcal{R}_{b}(X)$ and $\mathcal{R}_{c}(X)$. Using the
results of [7], one can obtain sufficient (but not necessary)
conditions for proximinality of these sets. This needs, besides
paths (see Section1), the consideration of some additional and
more complicated relations between points of $X$. The case $r\geq
3$ will not be considered in the current paper, since our main
purpose is to draw readers' attention to the arisen problems of
proximinality in the simplest case of approximation. For the
existing open problems connected with the set $\mathcal{R}\left(
\mathbf{a} ^{1},...,\mathbf{a}^{r}\right) $, where $r\geq 3$, see
[13] and [23].

\bigskip

\section*{1. Proximinality of $\mathcal{R}_{b}(X)$ in $B(X)$}

We begin this section with the definition of a path with respect
to two different directions $\mathbf{a}^{1}$ and $\mathbf{a}^{2}$.
A path with respect to the directions $\mathbf{a}^{1}$ and
$\mathbf{a}^{2}$ is a finite or infinite ordered set of points
$(\mathbf{x}^{1},\mathbf{x}^{2},...)$ in $ \mathbb{R}^{n}$ with
the units $\mathbf{x}^{i+1}-\mathbf{x}^{i}$, $i=1,2,..., $ in the
directions perpendicular alternatively to $\mathbf{a}^{1}$ and $
\mathbf{a}^{2}$. In the sequel, we simply use the term
\textquotedblleft path\textquotedblright\ instead of the long
expression \textquotedblleft path with respect to the directions
$\mathbf{a}^{1}$ and $\mathbf{a}^{2}$ \textquotedblright . The
length of a path is the number of its points and can be equal to
$\infty $ if the path is infinite. A singleton is a path of the
unit length. We say that a path $\left(
\mathbf{x}^{1},...,\mathbf{x} ^{m}\right) $ belonging to some
subset $X$ of $\mathbb{R}^{n}$ is irreducible if there is not
another path $\left( \mathbf{y}^{1},...,\mathbf{y }^{l}\right)
\subset X$ with $\mathbf{y}^{1}=\mathbf{x}^{1},~\mathbf{y}^{l}=
\mathbf{x}^{m}$ and $l<m$. If in a path $\left(
\mathbf{x}^{1},...,\mathbf{x} ^{m}\right) $ $m$ is an even number
and the set $\left( \mathbf{x}^{1},...,
\mathbf{x}^{m},\mathbf{x}^{1}\right) $ is also a path, then the
path $\left( \mathbf{x}^{1},...,\mathbf{x}^{m}\right) $ is called
to be closed. The notion of a path in the case when the directions
$\mathbf{a}^{1}$ and $ \mathbf{a}^{2}$ are basis vectors in
$\mathbb{R}^{2}$ was first introduced by Diliberto and Straus [6]
and exploited further in a number of works devoted to the
approximation of bivariate functions by univariate functions (see,
for example, [1,8,10,14,21]). Braess and Pinkus [2] used the
notion in their solution to one problem of interpolation by ridge
functions. It also appeared in problems of representation and well
approximation of a continuous multivariate function by functions
in $\mathcal{R}_{c}(X)$ (see [13]).

The following theorem follows from theorem 2 of [7]:

\smallskip

\textbf{Theorem 1.1.} \textit{Let $X\subset $ $\mathbb{R}^{n}$ and
the lengths of all irreducible paths in $X$ be uniformly bounded
by some positive integer. Then each function in $B(X)$ has an
extremal element in $ \mathcal{R}_{b}(X)$.}

\smallskip

There are a large number of sets in $\mathbb{R}^{n}$ satisfying
the hypothesis of this theorem. For example, if a set $X$ has a
cross section according to one of the directions $\mathbf{a}^{1}$
or $\mathbf{a}^{2}$, then the set $X$ satisfies the hypothesis of
theorem 1.1. By a cross section according to the direction
$\mathbf{a}^{1}$ we mean any set $X_{\mathbf{a} ^{1}}=\{x\in X:\
\mathbf{a}^{1}\cdot \mathbf{x}=c\},c\in \mathbb{R}$, with the
property: for any $\mathbf{y}\in X$ there exists a point
$\mathbf{y} ^{1}\in X_{\mathbf{a}^{1}}$ such that
$\mathbf{a}^{2}\cdot \mathbf{y}= \mathbf{a}^{2}\cdot
\mathbf{y}^{1}$. By the similar way, one can define a cross
section according to the direction $\mathbf{a}^{2}$. Regarding
theorem 1.1, one may ask if the condition of the theorem is
necessary for proximinality of $\mathcal{R}_{b}(X)$ in $B(X)$.
While we do not know a complete answer to this question, we are
going to give an example of a set $X $ for which theorem 1.1
fails. Let $\mathbf{a}^{1}=(1;-1),\ \mathbf{a} ^{2}=(1;1).$
Consider the set
\begin{eqnarray*}
X
&=&\{(2;\frac{2}{3}),(\frac{2}{3};-\frac{2}{3}),(0;0),(1;1),(1+\frac{1}{2}
;1-\frac{1}{2}),(1+\frac{1}{2}+\frac{1}{4};1-\frac{1}{2}+\frac{1}{4}),
\\
&&(1+\frac{1}{2}+\frac{1}{4}+\frac{1}{8};1-\frac{1}{2}+\frac{1}{4}-\frac{1}{8
}),...\}.
\end{eqnarray*}
In what follows, the elements of $X$ in the given order will be
denoted by $ \mathbf{x}^{0},\mathbf{x}^{1},\mathbf{x}^{2},...$ .
It is clear that $X$ is a path of the infinite length and
$\mathbf{x}^{n}\rightarrow \mathbf{x}^{0}$ , as $n\rightarrow
\infty $. Let $\sum_{n=1}^{\infty }c_{n}$ be any divergent series
with the terms $c_{n}>0$ and $c_{n}\rightarrow 0,$ as $
n\rightarrow \infty $. Besides let $f_{0}$ be a function vanishing
at the points $\mathbf{x}^{0},\mathbf{x}^{2},\mathbf{x}^{4},...,$
and taking values $c_{1},c_{2},c_{3},...$ at the points
$\mathbf{x}^{1},\mathbf{x}^{3},\mathbf{ x}^{5},...$ respectively.
It is obvious that $f_{0}$ is continuous on $X$. The set $X$ is
compact and satisfies all the conditions of proposition 2 of [21].
By this proposition, $\overline{\mathcal{R}_{c}(X)}=C(X).$
Therefore, for any continuous function on $X$, thus for $f_{0}$,
\begin{equation*}
\inf_{g\in \mathcal{R}_{c}(X)}\left\Vert f_{0}-g\right\Vert
_{C(X)}=0.\eqno (1.1)
\end{equation*}

Since $\mathcal{R}_{c}(X)\subset \mathcal{R}_{b}(X),$ we obtain
from (1.1) that

\begin{equation*}
\inf_{g\in \mathcal{R}_{b}(X)}\left\Vert f_{0}-g\right\Vert
_{B(X)}=0.\eqno (1.2)
\end{equation*}

Suppose that $f_{0}$ has an extremal element ${g_{1}^{0}\left(
\mathbf{a}^{1}{\cdot }\mathbf{x}\right) +g_{2}^{0}\left(
\mathbf{a}^{2}{ \cdot }\mathbf{x}\right) }$ in
$\mathcal{R}_{b}(X).$ By the definition of $ \mathcal{R}_{b}(X)$,
the ridge functions ${g_{i}^{0},i=1,2}$, are bounded on $X.$ From
(1.2)\ it follows that $f_{0}={g_{1}^{0}\left( \mathbf{a}^{1}{
\cdot }\mathbf{x}\right) +g_{2}^{0}\left( \mathbf{a}^{2}{\cdot
}\mathbf{x} \right) .}$ Since $\mathbf{a}^{1}\cdot
\mathbf{x}^{2n}=\mathbf{a}^{1}\cdot \mathbf{x}^{2n+1}$ and
$\mathbf{a}^{2}\cdot \mathbf{x}^{2n+1}=\mathbf{a} ^{2}\cdot
\mathbf{x}^{2n+2},$ for $n=0,1,...,$ we can write
\begin{equation*}
\sum_{n=0}^{k}c_{n+1}=\sum_{n=0}^{k}\left[
f(\mathbf{x}^{2n+1})-f(\mathbf{x} ^{2n})\right]
\end{equation*}

\begin{equation*}
=\sum_{n=0}^{k}\left[
{g_{2}^{0}}(\mathbf{x}^{2n+1})-{g_{2}^{0}}(\mathbf{x}
^{2n})\right] ={g_{2}^{0}(}\mathbf{a}^{2}\cdot
\mathbf{x}^{2k+1})-{g_{2}^{0}( }\mathbf{a}^{2}\cdot
\mathbf{x}^{0}).\eqno(1.3)
\end{equation*}

Since $\sum_{n=1}^{\infty }c_{n}=\infty ,$ we deduce from (1.3)\
that the function ${g_{2}^{0}\left( \mathbf{a}^{2}{\cdot
}\mathbf{x}\right) } $ is not bounded on $X.$ This contradiction
means that the function $f_{0}$ does not have an extremal element
in $\mathcal{R}_{b}(X).$ Therefore, the space $\mathcal{R}_{b}(X)$
is not proximinal in $B(X).$

\smallskip

\textit{Remark.} The above example is a slight generalization and
an adaptation of Havinson's example (see [10]) to our case.

\bigskip

\section*{2. Proximinality of $\mathcal{R}_{c}(X)$ in $C(X)$}

In this section, we are going to give sufficient conditions and
also a necessary condition for proximinality of
$\mathcal{R}_{c}(X)$ in $C(X).$

\smallskip

\textbf{Theorem 2.1.} \textit{Let the system of independent
vectors $ \mathbf{a}^{1}$ and $\mathbf{a}^{2}$ has a complement to
a basis $\{\mathbf{a }^{1},...,\mathbf{a}^{n}\}$ in
$\mathbb{R}^{n}$ with the property: for any point
$\mathbf{x}^{0}\in X$ and any positive real number $\delta $ there
exist a number $\delta _{0}\in (0,\delta ]$ and a point
$\mathbf{x}^{\sigma } $ in the set}

\begin{equation*}
\sigma =\{\mathbf{x}\in X:\mathbf{a}^{2}\cdot
\mathbf{x}^{0}-\delta _{0}\leq \mathbf{a}^{2}\cdot \mathbf{x}\leq
\mathbf{a}^{2}\cdot \mathbf{x}^{0}+\delta _{0}\},
\end{equation*}
\textit{such that the system}

\begin{equation*}
\left\{
\begin{array}{c}
\mathbf{a}^{2}\cdot \mathbf{x}^{\prime }=\mathbf{a}^{2}\cdot
\mathbf{x} ^{\sigma } \\ \mathbf{a}^{1}\cdot \mathbf{x}^{\prime
}=\mathbf{a}^{1}\cdot \mathbf{x} \\ \sum_{i=3}^{n}\left\vert
\mathbf{a}^{i}\cdot \mathbf{x}^{\prime }-\mathbf{a} ^{i}\cdot
\mathbf{x}\right\vert <\delta
\end{array}
\right. \eqno(2.1)
\end{equation*}
\textit{has a solution $\mathbf{x}^{\prime }\in \sigma $ for all
points $ \mathbf{x}\in \sigma .$Then the space
$\mathcal{R}_{c}(X)$ is proximinal in $ C(X).$}

\begin{proof} Introduce the following mappings and sets:
\begin{equation*}
\pi _{i}:X\rightarrow \mathbb{R}\text{, }\pi
_{i}(\mathbf{x)=a}^{i}\cdot \mathbf{x}\text{, }Y_{i}=\pi
_{i}(X\mathbf{)}\text{, }i=1,...,n.
\end{equation*}

Since the system of vectors
$\{\mathbf{a}^{1},...,\mathbf{a}^{n}\}$ is linearly independent,
the mapping $\pi =(\pi _{1},...\pi _{n})$ is an injection from $X$
into the Cartesian product $Y_{1} \times ...\times Y_{n} $ .
Besides, $\pi $ is linear and continuous. By the open mapping
theorem, the inverse mapping $\pi ^{-1}$ is continuous from $Y=\pi
(X)$ onto $X.$ Let $f$ be a continuous function on $X$. Then the
composition $f\circ \pi ^{-1}(y_{1},...y_{n})$ will be continuous
on $Y,$ where $y_{i}=\pi _{i}( \mathbf{x),}\ i=1,...,n,$ are the
coordinate functions. Consider the approximation of the function
$f\circ \pi ^{-1}$ by elements from
\begin{equation*}
G_{0}=\{g_{1}(y_{1})+g_{2}(y_{2}):\ g_{i}\in C(Y_{i}),\ i=1,2\}
\end{equation*}
over the compact set $Y$. Then one may observe that the function
$f$ has an extremal element in $\mathcal{R}_{c}(X)$ if and only if
the function $f\circ \pi ^{-1}$ has an extremal element in
$G_{0}$. Thus the problem of proximinality of $\mathcal{R}_{c}(X)$
in $C(X)$ is reduced to the problem of proximinality of $G_{0}$ in
$C(Y).$

Let $T,T_{1},...,T_{m+1}$ be metric compact spaces and $T\subset $
$ T_{1}\times ...\times T_{m+1}.$ For $i=1,...,m,$ let $\varphi
_{i}$ be the continuous mappings from $T$ onto $T_{i}.$ In [7],
the authors obtained sufficient conditions for proximinality of
the set

\begin{equation*}
C_{0}=\{\sum_{i=1}^{n}g_{i}\circ \varphi _{i}:\ g_{i}\in
C(T_{i}),\ i=1,...m\}
\end{equation*}
in the space $C(T)$ of continuous functions on $T.$ Since
$Y\subset $ $ Y_{1}\times Y_{2}\times Z_{3},$ where
$Z_{3}=Y_{3}\times ...\times Y_{n},$ we can use this result in our
case for the approximation of the function $ f\circ \pi ^{-1}$ by
elements from $G_{0}$. By this theorem, the set $G_{0}$ .is
proximinal in $C(Y)$ if for any $y_{2}^{0}\in Y_{2}$ and $\delta
>0$
there exists a number $\delta _{0}\in (0,$ $\delta )$ such that
the set $ \sigma (y_{2}^{0},\delta _{0})=[y_{2}^{0}-\delta
_{0},y_{2}^{0}+\delta _{0}]\cap Y_{2}$ has $(2,\delta )$ maximal
cross section. The last means that there exists a point
$y_{2}^{\sigma }\in \sigma (y_{2}^{0},\delta _{0})$ with the
property: for any point $(y_{1},y_{2},z_{3})\in Y,$ with the
second coordinate $y_{2}$ from the set $\sigma (y_{2}^{0},\delta
_{0}),$ there exists a point $(y_{1}^{\prime },y_{2}^{\sigma
},z_{3}^{\prime })\in Y$ such that $y_{1}=y_{1}^{\prime }$ and
$\rho (z_{3},z_{3}^{\prime })<\delta ,$ where $\rho $ is a metrics
in $Z_{3}.$ Since these conditions are equivalent to the
conditions of theorem 2.1, the space $G_{0}$ is proximinal in the
space $C(Y).$ Then by the above conclusion, the space
$\mathcal{R}_{c}(X)$ is proximinal in $C(X).$
\end{proof}

Let us give some simple examples of compact sets satisfying the
hypothesis of theorem 2.1. For the sake of brevity, we restrict
ourselves to the case $n=3.$

\begin{enumerate}
\item[(a)] Let $X$ be a closed ball in $\mathbb{R}^{3}$, $a^{1}$ and
$a^{2}$ be two arbitrary orthogonal directions. Then theorem 2.1
holds. Note that in this case, we can take $\delta _{0}=\delta $
and $a^{3}$ as an orthogonal vector to both the vectors $a^{1}$
and $a^{2}.$

\item[(b)] Let $X$ be the unite cube, $a^{1}=(1;1;0),\
a^{2}=(1;-1;0).$ Then theorem 2.1 also holds. In this case, we can
take $ \delta _{0}=\delta $ and $a^{3}=(0;0;1).$ Note that the
unit cube does not satisfy the hypothesis of the theorem for many
directions (take, for example, $a^{1}=(1;2;0)$ and$\
a^{2}=(2;-1;0)$).
\end{enumerate}

In the following example, one can not always chose $\delta _{0}$
as equal to $\delta $.

\begin{enumerate}
\item[(c)] Let $X=\{(x_{1},x_{2},x_{3}):\ (x_{1},x_{2})\in Q,\ 0\leq
x_{3}\leq 1\},$ where $Q$ is the union of two triangles
$A_{1}B_{1}C_{1}$ and $A_{2}B_{2}C_{2}$ with the vertices
$A_{1}=(0;0),\ B_{1}=(1;2),\ C_{1}=(2;0),\
A_{2}=(1\frac{1}{2};1),\ B_{2}=(2\frac{1}{2};-1),\ C_{2}=(3
\frac{1}{2};1).$ Let $a^{1}=(0;1;0)$ and $a^{2}=(1;0;0).$ Then it
is easy to see that theorem 2.1 holds (the vector $a^{3}$ can be
chosen as $(0;0;1)$). In this case, $\delta _{0}$ can not be
always chosen as equal to $\delta $. Take, for example,
$\mathbf{x}^{0}=(1\frac{3}{4};0;0)$ and $\delta =1\frac{3 }{4}.$
If $\delta _{0}=\delta ,$ then the second equation of the system
(2.1) has not a solution for a point $(1;2;0)$ or a point
$(2\frac{1}{2} ;-1;0).$ But if we take $\delta _{0}$ not more than
$\frac{1}{4}$, then for $ \mathbf{x}^{\sigma }=\mathbf{x}^{0}$ the
system has a solution. Note that the last inequality $\left\vert
\mathbf{a}^{3}\cdot \mathbf{x}^{\prime }- \mathbf{a}^{3}\cdot
\mathbf{x}\right\vert <\delta $ of the system can be satisfied
with the equality $\mathbf{a}^{3}\cdot \mathbf{x}^{\prime }=
\mathbf{a}^{3}\cdot \mathbf{x}$ if $a^{3}=(0;0;1).$
\end{enumerate}

It should be remarked that the results of [7] tell nothing about
necessary conditions for proximinality of the spaces considered
there. To fill this gap in our case, we want to give a necessary
condition for proximinality of $\mathcal{R}_{c}(X)$ in $C(X).$ Our
result will be based on the result of Marshall and O'Farrell given
below. First, let us introduce some notation. By
$\mathcal{R}_{c}^{i},\ i=1,2,$ we will denote the set of
continuous ridge functions $g\left( \mathbf{a}^{i}\cdot
\mathbf{x}\right) $ on the given compact set $X \subset
\mathbb{R}^{n}.$ Note that $
\mathcal{R}_{c}=\mathcal{R}_{c}^{1}+\mathcal{R}_{c}^{2}.$ Besides,
let $ \mathcal{R}_{c}^{3}=\mathcal{R}_{c}^{1}\cap
\mathcal{R}_{c}^{2}.$ For $ i=1,2,3,$ let $X_{i}$ be the quotient
space obtained by identifying points $ y_{1}$ and $y_{2}$ in $X$
whenever $f(y_{1})=f(y_{2})$ for each $f$ in $
\mathcal{R}_{c}^{i}.$ By $\pi _{i}$ denote the natural projection
of $X$ onto $X_{i},$ $i=1,2,3.$ Note that we have already dealt
with the quotient spaces $X_{1}$, $X_{2}$ and the projections $\pi
_{1},\pi _{2}$ in the previous section (see the proof of theorem
2.1). The relation on $X$, defined by setting $\ y_{1}\approx
y_{2}$ if $y_{1}$ and $y_{2}$ belong to some path, is an
equivalence relation. According to Marshall and O'Farrell [21] the
equivalence classes we call orbits. By $O(t)$ denote the orbit of
$X $ containing $t.$ For $Y\subset X,$ let $var_{Y}\ f$ be the
variation of a function $f$ on the set $Y.$ That is,

\begin{equation*}
\underset{Y}{var}f=\sup\limits_{x,y\in Y}\left\vert f\left(
x\right) -f\left( y\right) \right\vert .
\end{equation*}

\smallskip

\textbf{Theorem 2.2.} \textit{Suppose that the space $\mathcal{R}
_{c}(X)$ is proximinal in $C(X).$Then there exists a positive real
number c such that}

\begin{equation*}
\sup_{t\in X}\underset{O\left( t\right) }{var}\mathit{f\leq
c}\sup_{t\in X} \underset{\pi _{2}^{-1}\left( \pi _{2}\left(
t\right) \right) }{var}\mathit{f }\eqno(2.2)
\end{equation*}
\textit{for all $f$ in $\mathcal{R}_{c}^{1}.$}

\smallskip

The proof is simple. In [21], Marshall and O'Farrell proved the
following result (see Proposition 4 in [21]): Let $A_{1}\ $and
$A_{2}\ $be closed subalgebras of $C(X)\ $that contain the
constants. Let $(X_{1},\pi _{1}),\ (X_{2},\pi _{2})\ $and
$(X_{3},\pi _{3})\ $be the quotient spaces and projections
associated with the algebras $A_{1},$ $A_{2}\ $and $
A_{3}=A_{1}\cap A_{2}\ $respectively. Then $A_{1}+A_{2}\ $is
closed in $ C(X)\ $if and only if there exists a positive real
number $c$ such that
\begin{equation*}
\sup\limits_{z\in X_{3}}\underset{\pi _{3}^{-1}\left( z\right)
}{var}f\leq c\sup\limits_{y\in X_{2}}\underset{\pi _{2}^{-1}\left(
y\right) }{var}f\eqno (2.3)
\end{equation*}
for all $f\ $in $A_{1}.$

If $\mathcal{R}_{c}(X)$ is proximinal in $C(X),$ then it is
necessarily closed and therefore, by the above proposition, (2.3)
holds for the algebras $A_{1}^{i}=\mathcal{R}_{c}^{i},\ i=1,2,3.$
The right-hand side of (2.3) is equal to the right-hand side of
(2.2). Let $t$ be some point in $ X$ and $z=\pi _{3}(t).$ Since
each function $\ f\in \mathcal{R}_{c}^{3}$ is constant on the
orbit of $t$ (note that $f$ is both of the form ${ g_{1}\left(
\mathbf{a}^{1}{\cdot }\mathbf{x}\right) }$ and of the form ${
g_{2}\left( \mathbf{a}^{2}{\cdot }\mathbf{x}\right) }$),
$O(t)\subset \pi _{3}^{-1}(z).$ Hence,

\begin{equation*}
\sup_{t\in X}\underset{O\left( t\right) }{var}f\leq
c\sup\limits_{z\in X_{3}} \underset{\pi _{3}^{-1}\left( z\right)
}{var}f\eqno(2.4)
\end{equation*}

From (2.3) and (2.4) we obtain (2.2).

\smallskip

Note that the inequality (2.3) provides not worse but less
practicable necessary condition for proximinality than the
inequality (2.2) does. On the other hand, there are many cases in
which both the inequalities are equivalent. For example, let the
lengths of irreducible paths of $X$ are bounded by some positive
integer $n_{0}$. In this case, it can be shown that the inequality
(2.3), hence (2.2), holds with the constant $c=\frac{n_{0}}{2} $
and moreover $O(t)=\pi _{3}^{-1}(z)$ for all $t\in X$, where
$z=\pi _{3}(t) $ (see the proof of theorem 5 in [13]). Therefore,
the inequalities (2.2) and (2.3) are equivalent for the considered
class of sets $X.$ The last argument shows that all the compact
sets $X\subset $ $\mathbb{R}^{n}$ over which $\mathcal{R}_{c}(X)$
is not proximinal in $C(X)$ should be sought in the class of sets
having irreducible paths consisting sufficiently large number of
points. For example, let $I=[0;1]^{2}$ be the unit square, $
a^{1}=(1;1),\ a^{2}=(1;\frac{1}{2}).$ Consider the path

\begin{equation*}
l_{k}=\{(1;0),(0;1),(\frac{1}{2};0),(0;\frac{1}{2}),(\frac{1}{4};0),...,(0;
\frac{1}{2^{k}})\}.
\end{equation*}

It is clear that $l_{k}$ is an irreducible path with the length $
2k+2 $, where $k$ may be very large. Let $g_{k}$ be a continuous
univariate function on $\mathbb{R}$ satisfying the conditions:
$g_{k}(\frac{1}{2^{k-i}} )=i,\ i=0,...,k,$ $g_{k}(t)=0$ if
$t<\frac{1}{2^{k}},\ i-1\leq g_{k}(t)\leq i $ if $t\in
(\frac{1}{2^{k-i+1}},\frac{1}{2^{k-i}}),\ i=1,...,k,$ and $
g_{k}(t)=k$ if $t>1.$ Then it can be easily verified that

\begin{equation*}
\sup_{t\in X}\underset{\pi _{2}^{-1}\left( \pi _{2}\left( t\right)
\right) }{ var}g_{k}(\mathbf{a}^{1}{\cdot }\mathbf{x})\leq
1.\eqno(2.5)
\end{equation*}

Since $\max_{\mathbf{x}\in I}g_{k}(\mathbf{a}^{1}{\cdot
}\mathbf{x} )=k,$ $\min_{\mathbf{x}\in
I}g_{k}(\mathbf{a}^{1}{\cdot }\mathbf{x})=0$ and $
var_{\mathbf{x}\in O\left( t_{1}\right)
}g_{k}(\mathbf{a}^{1}{\cdot }\mathbf{ x})=k $ for $t_{1}=(1;0),$
we obtain that

\begin{equation*}
\sup_{t\in X}\underset{O\left( t\right)
}{var}g_{k}(\mathbf{a}^{1}{\cdot } \mathbf{x})=k.\eqno(2.6)
\end{equation*}

Since $k$ may be very large, from (2.5) and (2.6) it follows that
the inequality (2.2) cannot hold for the function
$g_{k}(\mathbf{a}^{1}{ \cdot }\mathbf{x})\in \mathcal{R}_{c}^{1}.$
Thus the space $\mathcal{R} _{c}(I)$ with the directions
$a^{1}=(1;1)$ and$\ a^{2}=(1;\frac{1}{2})$ is not proximinal in
$C(I).$

\smallskip

It should be remarked that if a compact set $X\subset $
$\mathbb{R} ^{n}$ satisfies the hypothesis of theorem 2.1, then
the length of all irreducible paths are uniformly bounded (see the
proof of theorem 2.1 and lemma in [7]). We have already seen that
if the last condition does not hold, then the proximinality of
both $\mathcal{R}_{c}(X)$ in $C(X)$ and $ \mathcal{R}_{b}(X)$ in
$B(X)$ fail for some sets $X.$ Besides the examples given above
and in Section 1, one can easily construct many other examples of
such sets. All these examples, theorems 1.1, 2.1, 2.2 and the
following remarks justify the statement of the following
conjecture:

\smallskip

\textbf{Conjecture. }\textit{Let $X$ be some subset of $\mathbb{R}
^{n}.$ The space $\mathcal{R}_{b}(X)$ is proximinal in $B(X)$ and
the space $ \mathcal{R}_{c}(X)$ is proximinal in $C(X)$ (in this
case, $X$ is considered to be compact) if and only if the lengths
of all irreducible paths of $X$ are uniformly bounded.}

\smallskip

\textit{Remark 1.} After completion of this work, Medvedev's
result came to our attention (see [16, p.58]). His result, in
particular, states that the set $R_{c}(X)$ is closed in $C(X)$ if
and only if the lengths of all irreducible paths of $X$ are
uniformly bounded. Thus, in the case of $ C(X)$, the necessity of
the above conjecture was proved by Medvedev.

\smallskip

\textit{Remark 2.} Note that there are situations in which a
continuous function (a specially chosen function on a specially
constructed set) has an extremal element in $\mathcal{R}_{b}(X)$,
but not in $\mathcal{R}_{c}(X)$ (see [16, p.73]). One subsection
of [16] (see p.68) devoted to the proximinality of sums of two
univariate functions with continuous and bounded summands in the
spaces of continuous and bounded bivariate functions respectively.
If $ X\subset \mathbb{R}^{2}$ and $\mathbf{a}^{1},\mathbf{a}^{2}$
be linearly independent directions in $\mathbb{R}^{2}$, then the
linear transformation $ y_{1}=$  $\mathbf{a}^{1}\cdot
\mathbf{x\,}$, $y_{2}=$ $\mathbf{a}^{2}\cdot \mathbf{x}$ reduces
the problems of proximinality of $\mathcal{R}_{b}(X)$ in $B(X)$
and $\mathcal{R}_{c}(X)$ in $C(X)$ to the problems considered in
that subsection. But in general, when $X\subset \mathbb{R}^{n},$
$n>2,$ our case cannot be obtained from that of [16].

\smallskip

\textbf{Acknowledgement.} I learned about the monograph by
Khavinson [16] from Allan Pinkus at the Technion. Using the
opportunity, I would like to express my sincere gratitude to him.

\end{document}